\documentclass[11pt,reqno]{amsart}
\usepackage[T1]{fontenc}
\usepackage{lmodern}
\usepackage{amsmath,amssymb,mathtools}
\usepackage[expansion=false]{microtype}
\usepackage[margin=30mm]{geometry}
\usepackage[hidelinks]{hyperref}
\hypersetup{pdftitle={Failure of endpoint estimates for lacunary bilinear spherical maximal functions on Heisenberg groups},
 pdfauthor={Sanghyuk Lee, Ji Li, Chaojie Wen, Liangchuan Wu},
 pdfkeywords={Heisenberg group, Hardy space, bilinear spherical means, lacunary maximal operator, endpoint counterexample}}
\allowdisplaybreaks
\newtheorem{theorem}{Theorem}[section]
\newtheorem{proposition}[theorem]{Proposition}
\newtheorem{lemma}[theorem]{Lemma}
\newtheorem{corollary}[theorem]{Corollary}

\theoremstyle{remark}
\newtheorem{remark}[theorem]{Remark}
\numberwithin{equation}{section}
\newcommand{\R}{\mathbb R}
\newcommand{\Hh}{\mathbb H}
\newcommand{\Hn}{{\mathbb H^n}}
\newcommand{\Sph}{\mathbb S}
\newcommand{\A}{\mathcal A}
\newcommand{\B}{\mathcal B}

\newcommand{\supp}{\operatorname{supp}}

\newcommand{\om}{\omega_0}
\newcommand{\Sh}{\mathfrak S}
\newcommand{\Mh}{\mathfrak M}
\title[The bilinear Hardy endpoint on the Heisenberg group]
{Failure of endpoint estimates for lacunary bilinear
spherical maximal functions on Heisenberg groups}
\author{Sanghyuk Lee, Ji Li, Chaojie Wen, Liangchuan Wu}
\makeatletter
\@ifundefined{subjclassname@2020}{
 \@namedef{subjclassname@2020}{\textup{2020} Mathematics Subject Classification}
}{}
\makeatother
\date{}
\begin{document}
\begin{abstract}
We prove that the lacunary bilinear Nevo--Thangavelu spherical
maximal operator on the Heisenberg group $\Hh^n$ does not map
$H^1(\Hh^n)\times H^1(\Hh^n)$ into
$L^{1/2,\infty}(\Hh^n)$ for any $n\geq1$.
We obtain this result by lifting the Euclidean counterexample in
\cite{EuclideanPaper} to the Heisenberg group.
\end{abstract}
\maketitle

\section{Introduction and the main result}
\label{hg:introduction}

In our recent paper \cite{EuclideanPaper}, we proved that the
lacunary bilinear spherical maximal operator on $\mathbb R^d$
does not map $H^1(\mathbb R^d)\times H^1(\mathbb R^d)$ into $L^{1/2,\infty}(\mathbb R^d)$ for any
$d\geq1$. In this note, we use the Euclidean construction to show
that the same endpoint estimate fails for the lacunary bilinear
horizontal spherical maximal operator on the Heisenberg group (studied by
Ghosh and Singh \cite{GhoshSingh}).
We give the details of the construction on the group.

Let $n\geq1$. We identify $\mathbb C^n$ with $\R^{2n}$ and put
$
 \om(z,w)=\operatorname{Im}(z\cdot\overline w).
$
The Heisenberg group $\Hh^n=\mathbb C^n\times\R$ has group law
(see, for example, \cite{FollandStein})
\begin{equation}\label{hg:group-law}
 (z,t)(w,s)=\bigl(z+w,t+s+\tfrac12\om(z,w)\bigr).
\end{equation}
Its Haar measure is $dz\,dt$, and its dilations are
$\delta_R(z,t)=(Rz,R^2t)$.
Let $d\sigma$ be normalized surface measure on
$\Sph^{4n-1}\subset\mathbb C^n\times\mathbb C^n$.
For $F,G\in C_c^\infty(\Hh^n)$, define
\begin{align}
 \Sh_R(F,G)(z,t)
 &=\int_{\Sph^{4n-1}}
 F\bigl((z,t)(Ry,0)^{-1}\bigr)
 G\bigl((z,t)(Rw,0)^{-1}\bigr)\,d\sigma(y,w)
 \label{hg:average}\\
 &=\int_{\Sph^{4n-1}}
 F\bigl(z-Ry,t-\tfrac R2\om(z,y)\bigr)
 G\bigl(z-Rw,t-\tfrac R2\om(z,w)\bigr)\,d\sigma(y,w).
 \nonumber
\end{align}
The same integral defines $\Sh_R$ for bounded, compactly supported
Borel functions. We use such representatives for the atoms below.

These are the bilinear horizontal spherical means introduced by
Ghosh and Singh \cite{GhoshSingh}.
The sphere in \eqref{hg:average} lies in the horizontal variables
of $\Hh^n\times\Hh^n$.
The central coordinates of the two translating elements are zero.
The terms involving $\om$ come from the group law.

For a finite set $E\subset(0,\infty)$, put
\begin{equation}\label{hg:maximal}
 \Mh_E(F,G)=\sup_{R\in E}|\Sh_R(F,G)|,\qquad
 \Mh_{\mathrm{lac}}(F,G)
 =\sup_{\ell\in\mathbb Z}|\Sh_{2^\ell}(F,G)|.
\end{equation}
Ghosh and Singh \cite[Theorem~1.5]{GhoshSingh} obtained bounds for
\eqref{hg:maximal} in an open range of Lebesgue exponents when
$n\geq2$. Their range does not include $(p_1,p_2)=(1,1)$.
At that point, the H\"older relation gives $p=1/2$.
We consider the weak estimate
\begin{equation}\label{hg:proposed}
 \|\Mh_{\mathrm{lac}}(F,G)\|_{L^{1/2,\infty}(\Hh^n)}
 \leq C\|F\|_{H^1(\Hh^n)}\|G\|_{H^1(\Hh^n)}.
\end{equation}
Here $H^1(\Hh^n)$ is the real Hardy space associated with the group
dilations; see \cite{FollandStein}. We use its atomic norm, as defined
in Section~\ref{hg:lifting}.
Recall also the weak quasinorm
\begin{equation}\label{eq:weak-quasinorm}
 \|h\|_{L^{p,\infty}(\Hh^n)}
 =\sup_{\lambda>0}\lambda
   |\{x\in \Hh^n:|h(x)|>\lambda\}|^{1/p}.
\end{equation}
In particular, the measure of a level set is squared when $p=1/2$.

\begin{theorem}\label{hg:main}
Let $n\geq1$. There is a constant $c_n>0$ with the following
property. For every sufficiently large integer $N$, there are a
set $E_N$ of $N$ dyadic radii and real functions
$F_N,G_N\in C_c^\infty(\Hh^n)$, each a finite linear combination
of smooth $(1,\infty)$-atoms in $H^1(\Hh^n)$, such that
\begin{equation}\label{hg:main-ratio}
 \frac{\|\Mh_{E_N}(F_N,G_N)\|_{L^{1/2,\infty}(\Hh^n)}}
 {\|F_N\|_{H^1(\Hh^n)}\|G_N\|_{H^1(\Hh^n)}}
 \geq c_n\log(N+1).
\end{equation}
Consequently, \eqref{hg:proposed} is false.
\end{theorem}

The dyadic radii in the proof form a geometric progression. Its
ratio is allowed to depend on $N$. The estimate
\eqref{hg:main-ratio} concerns this set of $N$ radii.
We do not assert that this logarithmic rate is sharp.

For $n=1$, our negative result uses the Euclidean construction
in dimension two and does not require a positive theorem on
$\Hh^1$.

\begin{corollary}\label{hg:L1}
For every $n\geq1$, $\Mh_{\mathrm{lac}}$ is not bounded from
$L^1(\Hh^n)\times L^1(\Hh^n)$ to
$L^{1/2,\infty}(\Hh^n)$.

\end{corollary}

{To prove Theorem~\ref{hg:main},} we use the Euclidean construction in our companion paper
\cite{EuclideanPaper}, with $d=2n$. Write
$$
 \A_R(f,g)(z)=\int_{\Sph^{4n-1}}
 f(z-Ry)g(z-Rw)\,d\sigma(y,w).
$$

In \cite{EuclideanPaper}, the construction begins with the
observation that the bilinear spherical average of two mean-zero
atoms need not cancel. To obtain a lower bound, we choose a fixed
function $\phi$ with $\int\phi=0$ and $\widehat\phi\geq0$.
The proof uses the ordinary Fourier transform on $\R^{2n}$ and,
in particular, the positivity of
$$
 \int_{\R}\widehat\phi(\tau y)\widehat\phi(\tau w)\,d\tau
$$
for suitable directions $y,w$. This leads to functions $f_N,g_N$
with small Hardy norms whose Euclidean maximal function is bounded
below on a set of positive measure; see
Proposition~\ref{hg:euclidean-witness}.

The Fourier transform on $\Hh^n$ is operator-valued, so the
Euclidean Fourier argument does not apply in the same form.
Instead, we lift the functions $f_N,g_N$ from the horizontal
space $\R^{2n}$ to $\Hh^n$.

The sphere in \eqref{hg:average} is
the same Euclidean sphere, but the group law introduces the terms
involving $\om$. To deal with them, for each atom centered at $u$
we use a cutoff in the coordinate
\begin{equation}\label{hg:q-coordinate}
 q_u(z,t)=t-\tfrac12\om(u,z).
\end{equation}
If $p=z-Ry$, then
$$
 q_u\bigl((z,t)(Ry,0)^{-1}\bigr)
 =t+\tfrac12\om(z-u,p).
$$
For $z$ in a fixed bounded set, $|t|\leq T/2$, and $T$ large
enough, the cutoff equals one whenever the atom contributes to the
average. The lifts $F_N,G_N$ then satisfy the exact identity
$$
 \Sh_R(F_N,G_N)(z,t)=T^{-2}\A_R(f_N,g_N)(z)
$$
{for such $(z,t)$; see Lemma~\ref{hg:equality-lemma}.}
It holds for every $R>0$ and includes all pairs of atoms and their
signs. Thus the Euclidean lower bound holds on a set of positive
measure in $\Hh^n$. {See \eqref{hg:final-equality} and \eqref{hg:final-level}.}

To estimate the Hardy norms, we split each lift into pieces with
central intervals of length comparable to the square of the atom
radius. Each piece is supported in a group ball and has integral
zero. The sum of their atomic coefficients is bounded independently
of the radius and of $T$. The two factors $T^{-1}$ in the average
are balanced by the length of the central interval in the weak
$L^{1/2}$ quasinorm.

Section~ \ref{hg:lifting} proves these two statements.
Section~\ref{hg:proof} proves Theorem~\ref{hg:main} with a central
cutoff fixed independently of $N$.

\section{Lifting from \texorpdfstring{$\R^{2n}$}{Euclidean space} to the Heisenberg group}
\label{hg:lifting}

Let $Q=2n+2$ denote the homogeneous dimension of $\Hh^n$.
For $x=(z,t)$, write
\begin{equation}\label{hg:gauge}
 B_\Hn(x_0,r)=\{x:|x_0^{-1}x|_\Hn<r\},
\end{equation}
where $|x|_\Hn=(|z|^4+t^2)^{1/4}$.
We use this homogeneous quasi-norm throughout. Replacing the
coefficient of $t^2$ by $16$ gives an equivalent gauge and the same
Hardy space with an equivalent norm.
The balls in \eqref{hg:gauge} are left translates of the ball at
the identity. Since $\delta_r$ has Jacobian $r^Q$, we have
\begin{equation}\label{hg:ball-volume}
 |B_\Hn(x_0,r)|=v_\Hn r^Q,
\end{equation}
where $v_\Hn=|B_\Hn (0,1)|$.
Write $B_E(u,r)$ for a Euclidean ball in $\R^{2n}$ and
let $v_E$ denote the volume of $B_E(0,1)$.

A $(1,\infty)$-atom in $H^1(\R^{2n})$ is a bounded function $a$ supported
in a ball $B_E(u,r)$, satisfying
\begin{equation}\label{hg:Euclidean-atom}
 \int_{\R^{2n}}a(z)\,dz=0,\qquad
 \|a\|_\infty\leq(v_Er^{2n})^{-1}.
\end{equation}
A $(1,\infty)$-atom in $H^1(\Hh^n)$ is a bounded function $b$ supported
in a ball $B_\Hn(x_0,r)$, satisfying
\begin{equation}\label{hg:Heisenberg-atom}
 \int_{\Hh^n}b(z,t)\,dz\,dt=0,\qquad
 \|b\|_\infty\leq(v_\Hn r^Q)^{-1}.
\end{equation}
The atomic Hardy space consists of the sums
$h=\sum_j\gamma_jb_j$, convergent in $L^1$, where each $b_j$ is a
$(1,\infty)$-atom in $H^1(\Hh^n)$ and $\sum_j|\gamma_j|<\infty$.
Its norm is the infimum of these coefficient sums.
This is the real Hardy space on the homogeneous group, with an
equivalent norm; see \cite[Theorem~3.30]{FollandStein}.
Only finite decompositions are used below.

Every atom in \eqref{hg:Heisenberg-atom} has $L^1$ norm at most
one. Therefore
\begin{equation}\label{hg:H1-L1}
 \|h\|_{L^1(\Hh^n)}\leq\|h\|_{H^1(\Hh^n)}.
\end{equation}
For a finite sum, each atomic decomposition gives an upper bound
for the Hardy norm. We shall estimate the coefficient sum of the
specific decomposition that we construct.

For $u\in\mathbb C^n$, set
\[
 q_u(z,t)=t-\tfrac12\om(u,z).
\]
Since $\om(u,u)=0$, the group law gives
\begin{equation}\label{hg:left-coordinates}
 (u,\tau)^{-1}(z,t)
 =\bigl(z-u,t-\tau-\tfrac12\om(u,z)\bigr)
 =\bigl(z-u,q_u(z,t)-\tau\bigr).
\end{equation}
For fixed $z$, the substitution $t\mapsto q_u(z,t)$ is a
translation. The full change of variables $(z,t)\mapsto(z,q_u)$
has Jacobian one.
In particular, $|z-u|\leq r$ and $|q_u(z,t)-\tau|\leq r^2$ imply
\begin{equation}\label{hg:cylinder-in-ball}
 |(u,\tau)^{-1}(z,t)|_\Hn\leq2^{1/4}r<2r.
\end{equation}
Thus the cylinder described above is contained in a group ball
of radius $2r$.

We now consider the lifting of $H^1(\R^{2n})$ atoms to $H^1(\Hn)$ functions. Fix a real function $\zeta\in C_c^\infty((-2,2))$ such that
\begin{equation}\label{hg:cutoff}
 0\leq\zeta\leq1,\qquad \zeta(v)=1\quad(|v|\leq1).
\end{equation}
For a $(1,\infty)$-atom $a$ in $H^1(\R^{2n})$ supported in $B_E(u,r)$ and $T>0$, define
\begin{equation}\label{hg:lift}
 \mathcal L_{u,T}a(z,t)
 =\frac1T a(z)\zeta\left(\frac{q_u(z,t)}T\right).
\end{equation}
The center $u$ is part of this definition.
When lifting a finite atomic sum, we apply \eqref{hg:lift}
separately to each atom.

\begin{lemma}\label{hg:lift-lemma}
Let $a$ satisfy \eqref{hg:Euclidean-atom}, and suppose $T\geq r^2$.
The lift \eqref{hg:lift} is a finite linear combination of
$(1,\infty)$-atoms in $H^1(\Hn)$, with coefficient sum at most $C_n$ independent of $u,r,T$, and $a$.
If $a$ is smooth, all atoms in this decomposition are smooth.
In particular,
\begin{equation}\label{hg:lift-norm}
 \|\mathcal L_{u,T}a\|_{H^1(\Hh^n)}\leq C_n.
\end{equation}
\end{lemma}

\begin{proof}
We split the coordinate $q_u$ into intervals of length comparable
to $r^2$. The cancellation of $a$ will hold in each piece.

Choose $\psi\in C_c^\infty((-1,1))$, $0\leq\psi\leq1$ such that
\[
 \sum_{k\in\mathbb Z}\psi(v-k)=1.
\]
Define
\begin{equation}\label{hg:pieces}
 A_k(z,t)=\frac1T a(z)\zeta\left(\frac{q_u(z,t)}T\right)
             \psi\left(\frac{q_u(z,t)}{r^2}-k\right).
\end{equation}
Then,
\begin{equation}\label{hg:finite-sum}
 \mathcal L_{u,T}a=\sum_{k\in\mathbb Z}A_k.
\end{equation}
If $A_k\ne0$, then $|q_u(z,t)|<2T$ and
$|q_u(z,t)-kr^2|<r^2$ at some point. Hence
\[
 |k|<\frac{2T}{r^2}+1.
\]
The integers satisfying this inequality lie in an interval of
length $4T/r^2+2$. Therefore
\begin{equation}\label{hg:number-pieces}
 \#\{k:A_k\ne0\}\leq\frac{4T}{r^2}+3.
\end{equation}
In particular, the sum in \eqref{hg:finite-sum} is finite.

We verify support, cancellation, and size for a fixed $k$.
On $\supp A_k$, $|z-u|\leq r$ and
$|q_u(z,t)-kr^2|\leq r^2$.
Equations \eqref{hg:left-coordinates}--\eqref{hg:cylinder-in-ball}
give
\begin{equation}\label{hg:piece-support}
 \supp A_k\subset B_\Hn((u,kr^2),2r).
\end{equation}

For cancellation, fix $z$. The map $t\mapsto q_u(z,t)=t-\tfrac12\om(u,z)$
is a translation, so
\begin{equation}\label{hg:vertical-integral}
 \int_{\R}\zeta\bigg(\frac{q_u(z,t)}{T}\bigg)\psi\bigg(\frac{q_u(z,t)}{r^2}-k\bigg)\,dt
 =\int_{\R}\zeta\bigg(\frac{q}{T}\bigg)\psi\bigg(\frac{q}{r^2}-k\bigg)\,dq,
\end{equation}
and the right-hand side does not depend on $z$. Both factors take
values in $[0,1]$, and $\zeta(q/T)=0$ for $|q|\geq2T$, so the
right-hand side is at most $4T$.

Since
$\int_{\Hh^n}|A_k(z,t)|\,dz\,dt\leq4\|a\|_{L^1(\R^{2n})}<\infty$,
Fubini's theorem and \eqref{hg:vertical-integral} give
\begin{equation}\label{hg:piece-cancellation}
 \int_{\Hh^n}A_k(z,t)\,dz\,dt
 =\frac1T\left(\int_{\R^{2n}}a(z)\,dz\right)
   \left(\int_{\R}\zeta\bigg(\frac{q}{T}\bigg)\psi\bigg(\frac{q}{r^2}-k\bigg)\,dq\right)=0.
\end{equation}
The first factor is zero, so every
piece has integral zero, including the pieces near the ends of the
cutoff.
Note also
\[
 \|A_k\|_\infty\leq\frac1{Tv_Er^{2n}}.
\]

Since $Q=2n+2$, we have
\begin{equation}\label{hg:piece-size}
 |B_\Hn((u,kr^2),2r)|\|A_k\|_\infty
 \leq\frac{v_\Hn 2^Q}{v_E}\frac{r^2}T.
\end{equation}

By \eqref{hg:piece-support}--\eqref{hg:piece-size}, each nonzero
$A_k$ is a $(1,\infty)$-atom multiplied by
$v_\Hn 2^Qr^2/(v_ET)$.
Thus \eqref{hg:finite-sum} is an atomic decomposition. By
\eqref{hg:number-pieces} and $r^2\leq T$, its coefficient sum is
at most

\begin{equation}\label{hg:coefficient-sum}
 \frac{v_\Hn 2^Q}{v_E}\,\frac{r^2}T\left(\frac{4T}{r^2}+3\right)
 \leq\frac{7\cdot2^Qv_\Hn }{v_E}
 =: C_n,
\end{equation}
which depends only on $n$.
Every factor in \eqref{hg:pieces} is smooth when $a$ is smooth.
This proves the lemma.
\end{proof}

The same substitution gives an exact $L^1$ identity:
\begin{equation}\label{hg:lift-L1}
 \|\mathcal L_{u,T}a\|_{L^1(\Hh^n)}
 =\|a\|_{L^1(\R^{2n})}\int_{\R}\zeta(v)\,dv.
\end{equation}
The Hardy estimate also uses the cancellation in
\eqref{hg:piece-cancellation} and the group balls in
\eqref{hg:piece-support}.

\begin{corollary}\label{hg:finite-lift}
Let $a_1,\dots,a_I$ and $b_1,\dots,b_J$ be $(1,\infty)$-atoms in $H^1(\R^{2n})$ as in
\eqref{hg:Euclidean-atom}, with $a_i$ supported in $B_E(u_i,r_i)$
and $b_j$ supported in $B_E(v_j,s_j)$, and let
$\alpha_1,\dots,\alpha_I$ and $\beta_1,\dots,\beta_J$ be real numbers.
Consider
\begin{equation}\label{hg:finite-horizontal}
 f=\sum_{i=1}^{I}\alpha_i a_i\qquad{\rm and}\qquad
 g=\sum_{j=1}^{J}\beta_j b_j.
\end{equation}
For $T\geq\max\{r_i^2,s_j^2:i,j\}$, put
\begin{equation}\label{hg:finite-vertical}
 F=\sum_{i=1}^{I}\alpha_i\mathcal L_{u_i,T}a_i\qquad{\rm and}\qquad
 G=\sum_{j=1}^{J}\beta_j\mathcal L_{v_j,T}b_j.
\end{equation}
These lifts depend on the decompositions in
\eqref{hg:finite-horizontal}, not only on $f$ and $g$.
Then
\begin{equation}\label{hg:finite-norm}
 \|F\|_{H^1(\Hh^n)}\leq C_n\sum_i|\alpha_i|\qquad{\rm and}\qquad
 \|G\|_{H^1(\Hh^n)}\leq C_n\sum_j|\beta_j|.
\end{equation}
If the atoms $a_i,b_j$ are smooth, $F,G$ are finite sums of
smooth, compactly supported $(1,\infty)$-atoms in $H^1(\Hn)$.
\end{corollary}

\begin{proof}
Apply Lemma~\ref{hg:lift-lemma} to every atom, multiply the
resulting coefficients by $\alpha_i$ or $\beta_j$, and add their
absolute values. There are finitely many atoms $a_i,b_j$, and
each decomposition in the lemma is finite.
This proves all assertions.
\end{proof}

We next choose $T$ so that the lifts recover the Euclidean
functions at every translated point that contributes to an average.
The Fourier calculation in \cite{EuclideanPaper} has already been
used to construct these Euclidean functions. The identity below is
proved directly from the group law.

For functions on $\R^{2n}$, define
\begin{equation}\label{hg:Euclidean-average}
 \B_E(f,g)=\sup_{R\in E}|\A_R(f,g)|,
\end{equation}
where $ \A_R(f,g)(z)
 =\int_{\Sph^{4n-1}}f(z-Ry)g(z-Rw)\,d\sigma(y,w)$.
This is the same sphere and the same normalized surface measure
as in \eqref{hg:average}.

The sign of the term $-\tfrac12\om(u,z)$ in $q_u$ is relevant. To check it,
write $p=z-Ry$. Since $Ry=z-p$ and $\om(z,z)=0$,
\begin{equation}\label{hg:translated-q}
 \begin{split}
 q_u\bigl((z,t)(Ry,0)^{-1}\bigr)
 &=t-\tfrac12\om(z,Ry)-\tfrac12\om(u,p)\\
 &=t+\tfrac12\om(z,p)-\tfrac12\om(u,p)\\
 &=t+\tfrac12\om(z-u,p).
 \end{split}
\end{equation}
Therefore
\begin{equation}\label{hg:single-translation}
 (\mathcal L_{u,T}a)\bigl((z,t)(Ry,0)^{-1}\bigr)
 =\frac1T a(p)
 \zeta\left(\frac{t+\frac12\om(z-u,p)}T\right).
\end{equation}

\begin{lemma}\label{hg:equality-lemma}
Consider the finite sums in \eqref{hg:finite-horizontal} and their
lifts in \eqref{hg:finite-vertical}. Suppose all centers have
Euclidean norm at most $U$, and all radii are at most $r_*$.
Let $M>0$, and choose
\begin{equation}\label{hg:T-choice}
 T\geq\max\{r_*^2,(M+U)(U+r_*)\}.
\end{equation}
For every $R>0$, $|z|\leq M$, and $|t|\leq T/2$,
\begin{equation}\label{hg:average-equality}
 \Sh_R(F,G)(z,t)=T^{-2}\A_R(f,g)(z).
\end{equation}
Consequently, for every finite set $E\subset(0,\infty)$,
\begin{equation}\label{hg:maximal-equality}
 \Mh_E(F,G)(z,t)=T^{-2}\B_E(f,g)(z)
 \qquad(|z|\leq M,\ |t|\leq T/2).
\end{equation}
\end{lemma}

\begin{proof}
Fix one atom $a_i$ and put $p=z-Ry$.
If $a_i(p)=0$, both sides of
\[
 (\mathcal L_{u_i,T}a_i)\bigl((z,t)(Ry,0)^{-1}\bigr)
 =T^{-1}a_i(z-Ry)
\]
are zero. If $a_i(p)\ne0$, its support gives $|p|\leq U+r_*$.
By $|\om(v,w)|\leq|v||w|$,
\begin{equation}\label{hg:shear-bound}
 \left|\tfrac12\om(z-u_i,p)\right|
 \leq\tfrac12(M+U)(U+r_*)\leq T/2.
\end{equation}
Together with $|t|\leq T/2$, this gives
$|t+\tfrac12\om(z-u_i,p)|\leq T$.
The cutoff in \eqref{hg:single-translation} is one by
\eqref{hg:cutoff}. Summing over $i$ gives
\begin{equation}\label{hg:full-F}
 F\bigl((z,t)(Ry,0)^{-1}\bigr)=T^{-1}f(z-Ry).
\end{equation}
For the second sum, put $p'=z-Rw$. On the support of $b_j(p')$,
$|p'|\leq U+r_*$ and
\[
 \left|\tfrac12\om(z-v_j,p')\right|
 \leq\tfrac12(M+U)(U+r_*)\leq T/2.
\]
Substituting in \eqref{hg:single-translation} for this atom and
summing over $j$ gives
\begin{equation}\label{hg:full-G}
 G\bigl((z,t)(Rw,0)^{-1}\bigr)=T^{-1}g(z-Rw).
\end{equation}
Multiplying \eqref{hg:full-F} and \eqref{hg:full-G}, then
integrating on $\Sph^{4n-1}$, proves \eqref{hg:average-equality}.
Taking absolute values and the supremum proves
\eqref{hg:maximal-equality}.
\end{proof}

\begin{remark}
The equality concerns the complete finite sums. In particular,
their signs are retained. Every pair of atoms occurs in
\eqref{hg:average-equality}. The proof does not discard pairs that
were not used to obtain the Euclidean lower bound.
No additional choice of signs is needed.
\end{remark}

\section{The counterexample on the Heisenberg group}
\label{hg:proof}

We use the following form of the Euclidean construction from
\cite[proof of Theorem~1.1, Sections~2--7]{EuclideanPaper}.
For $d\geq2$, we use $B_E$ for Euclidean balls in $\R^d$, and
$\A_R$ and $\B_E$ for the corresponding Euclidean averages and
maximal functions, defined as above with $2n$ replaced by $d$.

\begin{proposition}[\cite{EuclideanPaper}]\label{hg:euclidean-witness}
Let $d\geq2$. There are a real function
$\phi\in C_c^\infty(\R^d)$, supported in $\overline{B_E(0,2)}$
and satisfying $\int\phi=0$, and positive constants
$\rho,\lambda_*,b_*$ with the following property.
For every sufficiently large $N$, there are finite sets of centers
$u_i,v_j\in[0,1)^d$, a common $0<s<1$, positive numbers
$\mu_i,\nu_j$, and signs $\epsilon_i,\delta_j\in\{-1,1\}$ such that
\begin{gather}
 f_N(z)=\sum_i\epsilon_i\mu_i s^{-d}
                   \phi\left(\frac{z-u_i}s\right),\qquad
 g_N(z)=\sum_j\delta_j\nu_j s^{-d}
                   \phi\left(\frac{z-v_j}s\right),
 \label{hg:Euclidean-sums}\\
 \sum_i\mu_i=\sum_j\nu_j=\theta,\qquad
 \theta^2=\frac\rho{\log(N+1)}.
 \label{hg:Euclidean-coefficients}
\end{gather}
There are $N$ dyadic radii $E_N$ and a measurable set
$W_N\subset[-1,2]^d$ satisfying
\begin{equation}\label{hg:Euclidean-witness-level}
 |W_N|\geq b_*,\qquad
 \B_{E_N}(f_N,g_N)(z)>\lambda_*\quad(z\in W_N).
\end{equation}
All constants and the function $\phi$ are independent of $N$.
The radii can be taken as
\begin{equation}\label{hg:radii}
 E_N=\{\mathcal L^{j-N}:1\leq j\leq N\},\qquad
 \mathcal L\geq\max\{2^7\sqrt d,N^{1/d}\},
\end{equation}
where $\mathcal L$ is a dyadic integer.
Every $R\in E_N$ is at most one.
\end{proposition}

These are the finite functions constructed in
\cite[proof of Theorem~1.1]{EuclideanPaper}.
The radii $\mathcal L^{j-N}$ are defined in Section~2 of that paper, and the
coefficient identities \eqref{hg:Euclidean-coefficients} are proved in Section~3.
The sign selection $\epsilon_i,\delta_j$ in Section~7 gives the lower bound on a subset
of the union of the selected regions. The containing regions in
Section~6 lie in $[-1,2]^d$, so the set may be chosen there.

\subsection{Explicit lifts and their atomic decompositions}

We use Proposition~\ref{hg:euclidean-witness} with $d=2n$.
Define
\begin{align}
 F_N(z,t)
 &=\sum_i\frac{\epsilon_i\mu_i}{T}s^{-d}
   \phi\left(\frac{z-u_i}s\right)
   \zeta\left(\frac{t-\frac12\om(u_i,z)}T\right),
 \label{hg:explicit-F}\\
 G_N(z,t)
 &=\sum_j\frac{\delta_j\nu_j}{T}s^{-d}
   \phi\left(\frac{z-v_j}s\right)
   \zeta\left(\frac{t-\frac12\om(v_j,z)}T\right)
 \label{hg:explicit-G}
\end{align}
with $ T=100(d+1)$.
These are real smooth functions with compact support.
For example, every nonzero summand of $F_N$ satisfies
\[
 |z-u_i|\leq2s\leq2,\qquad |q_{u_i}(z,t)|\leq2T.
\]
Since $|u_i|\leq\sqrt d$, this gives
\begin{equation}\label{hg:fixed-support}
 |z|\leq\sqrt d+2,\qquad
 |t|\leq2T+\tfrac12\sqrt d(\sqrt d+2).
\end{equation}
The same bounds hold for $G_N$.
The set defined by these two non-strict inequalities is compact
and independent of $N$.

To apply Lemma~\ref{hg:lift-lemma}, choose a constant $C_\phi>0$
such that
$
 C_\phi\geq v_E3^d\|\phi\|_\infty.
$
Then
$$
 a_i(z)=C_\phi^{-1}s^{-d}\phi\left(\frac{z-u_i}{s}\right)
$$
is a smooth $(1,\infty)$-atom in $H^1(\R^d)$ supported in $B_E(u_i,3s)$, and the
same holds for
$$
 b_j(z)=C_\phi^{-1}s^{-d}\phi\left(\frac{z-v_j}{s}\right).
$$

The coefficients of these
atoms in $f_N,g_N$ are $C_\phi\epsilon_i\mu_i$ and
$C_\phi\delta_j\nu_j$, respectively. Since $(3s)^2\leq9<T$,
Lemma~\ref{hg:lift-lemma} and \eqref{hg:Euclidean-coefficients} give
\begin{equation}\label{hg:final-Hardy}
 \|F_N\|_{H^1(\Hh^n)},\|G_N\|_{H^1(\Hh^n)}
 \leq C_n\theta.
\end{equation}
Both functions are finite sums of smooth $(1,\infty)$-atoms in $H^1(\Hh^n)$.

\subsection{The level set and the norm comparison}

We check that the same fixed $T$ gives equality of the averages
on the required set. Since $W_N\subset[-1,2]^d$, its points
satisfy $|z|\leq2\sqrt d$.
The centers satisfy $|u_i|,|v_j|\leq\sqrt d$, and the supporting
balls above have radii at most three.
With $M=2\sqrt d$, $U=\sqrt d$, and $r_*=3$,
\begin{equation}\label{hg:fixed-T-check}
 (M+U)(U+r_*)=3\sqrt d(\sqrt d+3)\leq10d<T,\qquad
 r_*^2=9<T.
\end{equation}
Thus \eqref{hg:T-choice} holds.
Lemma~\ref{hg:equality-lemma} applies to the atomic sums above and
yields
\begin{equation}\label{hg:final-equality}
 \Sh_R(F_N,G_N)(z,t)=T^{-2}\A_R(f_N,g_N)(z)
 \quad(z\in W_N,\ |t|\leq T/2,\ R>0).
\end{equation}
Taking the supremum over the same $N$ radii and using
\eqref{hg:Euclidean-witness-level}, we obtain
\begin{equation}\label{hg:final-level}
 W_N\times[-T/2,T/2]
 \subset\{x:\Mh_{E_N}(F_N,G_N)(x)>\lambda_*/T^2\}.
\end{equation}
Since Haar measure is $dz\,dt$,
\begin{align*}
 \left|\left\{x\in\Hh^n:
 \Mh_{E_N}(F_N,G_N)(x)>\frac{\lambda_*}{T^2}
 \right\}\right|
 &\geq \left|W_N\times[-T/2,T/2]\right|
 =T|W_N|
 \geq Tb_*.
\end{align*}
Taking $\lambda=\lambda_*/T^2$ in
\eqref{eq:weak-quasinorm}, we obtain
\begin{align}
 \|\Mh_{E_N}(F_N,G_N)\|_{L^{1/2,\infty}}
 &\geq \frac{\lambda_*}{T^2}
 \left|\left\{x\in\Hh^n:
 \Mh_{E_N}(F_N,G_N)(x)>\frac{\lambda_*}{T^2}
 \right\}\right|^2 \notag\\
 &\geq\frac{\lambda_*}{T^2}(Tb_*)^2
 =\lambda_*b_*^2.
 \label{hg:final-weak}
\end{align}
In particular, $F_N$ and $G_N$ are nonzero. Together with
\eqref{hg:final-Hardy} and $\theta^2=\rho/\log(N+1)$, this gives
\begin{align*}
 \frac{\|\Mh_{E_N}(F_N,G_N)\|_{L^{1/2,\infty}}}
 {\|F_N\|_{H^1}\|G_N\|_{H^1}}
 &\geq\frac{\lambda_*b_*^2}{C_n^2\theta^2}
 =\frac{\lambda_*b_*^2}{C_n^2\rho}\log(N+1).
\end{align*}
Since $d=2n$ and the constants are independent of $N$, this proves
\eqref{hg:main-ratio} and completes the proof of
Theorem~\ref{hg:main}.

\begin{remark}

A similar approach may extend to other groups once suitable bilinear spherical maximal functions are defined in that setting. We do not pursue this here.

\end{remark}

\bigskip
\bigskip

\noindent \textbf{Acknowledgements.} S. Lee is supported by the National Research Foundation of Korea (NRF) through Grant No. RS-2024-00342160.  J. Li is supported by ARC DP260100485. Liangchuan Wu is supported by NNSF of China \#12201002.

\bigskip
\bigskip
\bigskip
\noindent (S. Lee) Department of Mathematical Sciences and RIM,
Seoul National University, Seoul 08826, Republic of Korea\par
\noindent\textit{E-mail}: \texttt{shklee@snu.ac.kr}

\bigskip
\noindent (J. Li) School of Mathematical and Physical Sciences,
Macquarie University, NSW, 2109, Australia\par
\noindent\textit{E-mail}: \texttt{ji.li@mq.edu.au}

\bigskip
\noindent (C. Wen) Department of Mathematics, Sun Yat-sen University,
Guangzhou, 510275, China\par
\noindent\textit{E-mail}: \texttt{wenchj@mail2.sysu.edu.cn}

\bigskip
\noindent (L. Wu) School of Mathematical Sciences, Anhui University,
Hefei, 230601, P.R.~China\par
\noindent\textit{E-mail}: \texttt{wuliangchuan@ahu.edu.cn}

\end{document}